\documentclass[12pt, reqno]{amsart}
\usepackage{amsthm,amsfonts,amssymb,euscript}

\theoremstyle{plain}
  \newtheorem{theorem}[subsection]{Theorem}

 \theoremstyle{remark}

\theoremstyle{definition}

\begin{document}

\title[Counterexamples to $L^p$ collapsing]%
{Counterexamples to $L^p$ collapsing estimates}

\author[X. Du]{Xiumin Du}
\address{University of Maryland, College Park}
\email{xdu@math.umd.edu}

\author[M. Machedon]{Matei Machedon}
\address{University of Maryland, College Park}
\email{mxm@math.umd.edu}

\date{\today}
\dedicatory{}
\commby{}

\maketitle
\begin{abstract}
We show that certain $L^2$ space-time estimates for generalized density matrices which have been used by several authors in recent years to study equations of BBGKY or Hartree-Fock type, do not have non-trivial $L^pL^q$ generalizations.
\end{abstract}
\section{Introduction and main results}

In recent years, effective equations approximating the evolution of a large number of interacting Bosons or Fermions have been studied extensively. The best known example is the celebrated work of Erd\"os, Schlein and Yau
\cite{E-S-Y2},
\cite{E-S-Y3}.

Since that work, a number of authors have studied the related
 Gross-Pitaevskii or BBGKY hierarchies, or the Hartree-Fock or Hartree-Fock-Bogoliubov equations, using harmonic analysis techniques and space-time
$L^2$ estimates for a suitable trace density of solutions of the linear Schr\"odinger equation. We call such estimates ``collapsing estimates'', and list several instances, all in 3 space dimensions (thus, $x \in \mathbb R^3$, etc.).

If
\begin{align}
G(t, x, y, z)=e^{ \frac{it (\Delta_x + \Delta_y-\Delta_z )}{2}}G_0, \label{G}
\end{align}
 then
\begin{align}
&\|\nabla_x G(t, x, x, x) \|_{L^2(dt dx)}
\lesssim \|\nabla_x \nabla_y \nabla_z G_0(x, y, z)\label{++-}
\|_{L^2(dx dy dz)}.
\end{align}
This estimate was used in the study of  the Gross-Pitaevskii or BBGKY hierarchies.
See \cite{K-MMM} (where the estimate originates), as well as \cite{C-P}, \cite{C-H1},
\cite{C-H2}.

Another related example is: if
\begin{align}
\Lambda(t, x, y)=e^{ \frac{it (\Delta_x + \Delta_y)}{2}}\Lambda_0, \label{Lambda}
\end{align}
then
\begin{align}
& \||\nabla|_x^{1/2} \Lambda(t, x, x) \|_{L^2(dt dx )} \lesssim \||\nabla|_x^{1/2} |\nabla|_y^{1/2} \Lambda_0(x, y)\|_{L^2(dx dy)}.
\label{++}
\end{align}
This estimate is useful for the Hartree-Fock-Bogoliubov equations, see
\cite{GM1}, \cite{GM2}.

Finally, if

\begin{align}
\Gamma(t, x, y)=e^{ \frac{it (\Delta_x - \Delta_y)}{2}}\Gamma_0, \label{Gamma}
\end{align}
then
\begin{align}
&\||\nabla_x|^{\frac{1}{2}} \langle \nabla_x\rangle^{2 \epsilon}\Gamma(t, x, x) \|_{L^2(dt dx)}
&\lesssim_{\epsilon} \|\langle\nabla_x\rangle^{\frac{1}{2}+\epsilon}\langle\nabla_y\rangle^{\frac{1}{2}+\epsilon}  \Gamma_0(x, y)\|_{L^2(dx dy)}. \label{+-}
\end{align}
Such estimates are relevant to both the Hartree-Fock-Bogoliubov equations mentioned above,
and Hartree-Fock. See Theorem 3.3 in \cite{C-H-P}.

We also mention the approach of \cite{FLLS}, \cite{FS} which applies to equation \eqref{Gamma} and allows a wide range of $L^p(dt) L^q(dx)$ estimates on the left hand side, but the right hand side of the inequality is estimated
in a Schatten norm.

It is natural to ask whether one can replace the $L^2(dt) L^2(dx)$ norm on the left hand side of  estimates \eqref{++-},
\eqref{++} or \eqref{+-}
by an $L^p(dt) L^q(dx)$ norm, while keeping the right hand side in a Sobolev norm, which is useful for applications to PDEs. One can trivially make $p$ or $q$ bigger than $2$ by putting more derivatives on the right hand side, so the interesting question is if one can make $p$ or $q$ less than 2.

The main result of this note is that this is impossible.

We prove the following closely related results.

\begin{theorem}\label{lambdatheorem}
Let $\Lambda$ be given by \eqref{Lambda}, with $x, y \in \mathbb R^n$. Assume
\begin{align}
& \||\nabla|_x^{\alpha} \Lambda(t, x, x) \|_{L^p(dt)L^q( dx )} \lesssim \| \Lambda_0(x, y)\|_{H^s(dx dy)}
\label{++pq}
\end{align}
for some $\alpha \ge 0, s \ge 0$. Then $p \ge 2$ and $q \ge 2$.
\end{theorem}
\begin{theorem}\label{gammatheorem}
Let $\Gamma$ be given by \eqref{Gamma}, with $x, y \in \mathbb R^n$. Assume
\begin{align}
& \||\nabla|_x^{\alpha} \Gamma(t, x, x) \|_{L^p(dt)L^q( dx )} \lesssim \| \Gamma_0(x, y)\|_{H^s(dx dy)}
\label{+-pq}
\end{align}
for some $\alpha \ge 0, s \ge 0$. Then $p \ge 2$ and $q \ge 2$.
\end{theorem}
\begin{theorem}\label{Gtheorem}
Let $G$ be given by \eqref{G}, with $x, y, z \in \mathbb R^n$. Assume
\begin{align}
& \||\nabla|_x^{\alpha} G(t, x, x, x) \|_{L^p(dt)L^q( dx )} \lesssim \| G_0(x, y, z)\|_{H^s(dx dy dz)}
\label{++-pq}
\end{align}
for some $\alpha \ge 0, s \ge 0$. Then $p \ge 2$ and $q \ge 2$.
\end{theorem}

\textbf{Acknowledgements.} The first author is supported by the National Science Foundation under Grant No. DMS-$1856475$.

\section{Proofs}
\subsection{Proof of Theorem \ref{lambdatheorem}}

\subsubsection{Necessity of $p\geq 2$.} \label{sec-lambdap}

Let $R$ be a large number (which will approach $\infty$ at the end of the proof). Let $C$ be a fixed large number (depending on $n$).
Let
\begin{align*}
F_0(x, y)=  e^{-\frac{|x|^2+|y|^2}{2 C R}}
\end{align*}
so that
\begin{align}
e^{ \frac{it (\Delta_x +\Delta_y)}{2}}F_0:=
F(t, x, y)= \frac{1}{(1+i t/(CR))^n} e^{-\frac{|x|^2+|y|^2}{2(CR+it)}}. \label{tube}
\end{align}
We think of
$F(t, x, y)$ as the basic ``vertical tube'' solution
to the linear Schr\"odinger equation in $2n+1$ dimensions
which is essentially 1 if $|x|, |y| \le R^{1/2}$, $0 \le t \le R$.
The rigorous statement is that $C$ is chosen so that $\Re F(t, x, y) \ge \frac{1}{2}$ in the above range. Also, the Fourier transform (in space) of $F$ is essentially supported at frequencies $|\xi|, |\eta| \le R^{-1/2}$.

We choose the function $\Lambda(t, x, y)$ to be a sum of translates and modulations of $F(t, x, y)$ which are inclined at 45 degrees and are trained to reach the region
$|x|\le \frac{1}{100}$, $|y| \le\frac{1}{100}$, $R-R^{\frac{1}{2}}<t<R$ with almost the same oscillation (and almost no cancellations). The summands will have Fourier transforms essentially supported in balls of radius $ R^{-1/2}$
centered at unit vectors.

 Explicitly,
 choose roughly $R^{n-\frac{1}{2}}$ points $(x_k, y_k)$ which are spaced at distance $R^{1/2}$ from each other
on the sphere  $|(x, y)|=R$. For technical reasons, we only choose points for which all coordinates are $\ge \frac{R}{10 n }$.
Define $(\xi_k, \eta_k)=\frac{(x_k, y_k)}{R}$.

Choose the following initial conditions:
$$
\Lambda_0(x, y)=\sum e^{i (x \cdot \xi_k + y \cdot \eta_k)}F_0(x+x_k, y+y_k).
$$
The functions being summed are approximately orthogonal and each have $L^2$ norm $\sim R^{n/2}$:
\begin{equation}
\int \big|F_0(x+x_k, y+y_k)F_0(x+x_l, y+y_l)\big| dxdy
= \pi^n (CR)^n e^{-\frac{|(x_k, y_k)-(x_l, y_l)|^2}{4 C R}}. \label{orthogonal}
\end{equation}
Recalling that the sum has $\sim R^{n-\frac{1}{2}}$ terms, we derive
\begin{align*}
\|\Lambda_0\|_{L^2(dx dy)} \lesssim  R^{n-\frac{1}{4}}.
\end{align*}
The same type of upper bound holds for higher order derivatives (since $|(\xi_k, \eta_k)|=1$), thus, for each fixed $s$,
\begin{align}
\|\Lambda_0\|_{H^s(dx dy)} \lesssim  R^{n-\frac{1}{4}}. \label{hsnorm}
\end{align}

The solution looks like
\begin{align*}
\Lambda(t, x, y)&=\sum e^{-i t \frac{ (|\xi_k|^2+|\eta_k|^2)}{2}}e^{i (x \cdot \xi_k + y \cdot \eta_k)}F(t, x+x_k- t \xi_k, y+y_k-t \eta_k)\\
&=e^{-i  \frac{t}{2}}\sum e^{i (x \cdot \xi_k + y \cdot \eta_k)}F(t, x+x_k- t \xi_k, y+y_k-t \eta_k),
\end{align*}
and
\begin{align*}
\big|\Lambda(t, x, y)\big|\ge \Re \sum e^{i (x \cdot \xi_k + y \cdot \eta_k)}F(t, x+x_k- t \xi_k, y+y_k-t \eta_k) \sim
R^{n-\frac{1}{2}},
\end{align*}
 if $|(x, y) | \le \frac{1}{100}$, $ R-R^{\frac{1}{2}}<t<R$.
Thus \begin{align}
  R^{\frac{1}{2p}} R^{n- \frac{1}{2}} \lesssim \|\Lambda(t, x, x)\|_{L^p(dt)L^q(dx)}, \label{mainest}
\end{align}
so, recalling \eqref{hsnorm}, if
\begin{align*}
& \| \Lambda(t, x, x) \|_{L^p(dt)L^q( dx )} \lesssim \| \Lambda_0(x, y)\|_{H^s(dx dy)},
\end{align*}
then $p \ge 2$.

Using the product rule and the lower bounds on the components of $\xi_k, \eta_k$, same argument works for ordinary derivatives of order $\alpha =m \in \mathbb N$.

To justify the statement for fractional derivatives of non-integer order $\alpha$, do
a Littlewood-Paley decomposition in space $\Lambda(t, \cdot, \cdot)=
P_{\le 10}\Lambda(t, \cdot, \cdot) + P_{\ge 10}\Lambda(t, \cdot, \cdot)$,
where $P_{\le 10}$ localizes functions of $2n $ variables, smoothly at frequencies $\le 10$.
Then $P_{\ge 10}\Lambda(t, \cdot, \cdot)$ is exponentially small as $R\to \infty$. This is true for the function $F_0$, and its translates by a unit vector in Fourier space.

 A crude estimate is
\begin{align*}
\|P_{\ge 10}\Lambda(t, \cdot, \cdot)\|_{H^s} \lesssim_s e^{- \sqrt R}.
\end{align*}
 For our counterexample, we use $P_{\le 10}\Lambda(t, \cdot, \cdot)$ instead of
 $\Lambda(t, \cdot, \cdot)$.

Thus, for $R$ sufficiently large, $
|\nabla^m P_{\le 10}\Lambda(t, x, y)| \sim |\nabla^m \Lambda(t, x, y)|
\sim R^{n-\frac{1}{2}}
$ if $|(x, y) | \le \frac{1}{100}$, $ R-R^{\frac{1}{2}}<t<R$.
The function $\left(P_{\le 10}\Lambda\right)(t, x, x)$ is
supported, in Fourier  space, at frequencies $|\xi| \le 20$. Denote, by abuse of notation, $P_{\le 20}$  the operator localizing functions of $n$ variables at frequencies $|\xi| \le 20$. Let $m \in \mathbb N$, $m > \alpha$. Then the operator $\frac{\nabla^m}{|\nabla|^{\alpha}}P_{\le 20}$ (defined in the obvious way on the Fourier transform side) is bounded on all $L^p$ spaces, and
 \begin{align*}
  R^{\frac{1}{2p}} R^{n- \frac{1}{2}}& \lesssim \| \nabla^m
 \left( P_{\le 10}\Lambda\right)(t, x, x)\|_{L^p(dt)L^q(dx)} \\
 &=\|\frac{\nabla^m}{|\nabla|^{\alpha}}P_{\le 20} |\nabla|^{\alpha}\left( P_{\le 10}\Lambda\right)(t, x, x)\|_{L^p(dt)L^q(dx)}\\
 &\lesssim \||\nabla|^{\alpha} \left( P_{\le 10}\Lambda\right)(t, x, x)\|_{L^p(dt)L^q(dx)},
\end{align*}
while
\begin{align*}
\|P_{\le 10}\Lambda_0\|_{H^s(dx dy)} \lesssim C^n R^{n-\frac{1}{4}}.
\end{align*}
Letting $R \to \infty$,  we conclude $p\ge 2$ as before.

\subsubsection{Necessity of $q \geq 2$.} \label{sec-lambdaq}
 Let $F(t, x, y)$ be the basic vertical tube solution of height $R$ (as \eqref{tube}). Let $m\gg1$. Choose roughly $R^{mn-\frac n2}$ points $x_k$ which are spaced at distance $\sim R^{\frac{1}{2}}$ in a large ball $B(0, R^m)$ of radius $R^m$ in $\mathbb R^n$. Fix a unit vector $\xi\in S^{n-1}$.

 We take initial conditions
 \begin{align*}
\Lambda_0( x, y)=e^{i (x+y)\cdot \xi}\sum F_0(x+x_k,y+x_k).
\end{align*}
Then
\begin{align*}
\Lambda(t, x, y)=e^{i (x+y)\cdot \xi}e^{- i t}\sum F(t, x+x_k-t\xi, y+x_k-t\xi).
\end{align*}

There are roughly $R^{mn - \frac{n}{2}}$ terms in the sum. The summands are essentially orthogonal (as in \eqref{orthogonal} ) and each term has $L^2$ norm
$\sim R^{n/2}$, thus
\begin{align*}
\|\Lambda_0\|_{L^2(dxdy)}\sim  R^{\frac{n}{4}+\frac{mn}{2}}.
\end{align*}
On the other hand, each $F(t, x+x_k-t\xi, y+x_k-t\xi)$ is essentially $1$ on a tube $T_k$ of radius $R^{1/2}$ and length $R$ in $2n+1$ dimensions, and rapidly decaying out of $T_k$. Note that at $t=0$, $T_k$ is centered at $(0, -x_k,-x_k)$. Moreover, these tubes $T_k$ are in the same direction $(1,\xi,\xi)$ and hence disjoint. Therefore,
$|\Lambda(t,x,y)|\gtrsim 1$
on the union of the tubes $T_k$. In particular, $|\Lambda(t,x,x)|\gtrsim 1$ for $0\leq t\leq R$ and $x\in B(t\xi, R^m)$.
 We only need the previous estimate for $0\leq t\leq 1$, where the claim is obvious.
In addition, the Fourier transform of $\Lambda (t, x, x)$ is supported (essentially) in a $R^{-\frac{1}{2}}$ neighbourhood of the point $2 \xi$ , with $|\xi|=1$,
so $||\nabla|^{\alpha}\Lambda(t,x,x)|\gtrsim 1$ for $0\leq t\leq 1$ and $x\in B(t\xi, R^m)$.
Thus
$$\||\nabla|^{\alpha}\Lambda(t, x, x)\|_{L^p([0, 1])L^q(dx)}\gtrsim  R^{\frac{mn}{q}},$$
while $\|\Lambda_0\|_{H^s(dxdy)}\sim
\|\Lambda_0\|_{L^2(dxdy)}\sim R^{\frac{n}{4}+\frac{mn}{2}}$ and $m\gg 1$, so $q \ge 2$ is necessary.

\subsection{Proof of Theorem \ref{gammatheorem}}

The examples for $\Gamma$ are similar to those for $\Lambda$, and are included for completeness.

\subsubsection{Necessity of $p\geq 2$.}

First we take the basic ``vertical tube'' solution. Let
\begin{align*}
F_0(x, y)=  e^{-\frac{|x|^2+|y|^2}{2 C R}}
\end{align*}
so that
\begin{align}
e^{ \frac{it (\Delta_x -\Delta_y)}{2}}F_0:=
F(t, x, y)= \frac{1}{(1+(\frac{t}{CR})^2)^{\frac n2}} e^{-\frac{|x|^2}{2(CR+it)}}e^{-\frac{|y|^2}{2(CR-it)}}. \label{gammatube}
\end{align}
The solution
$F(t, x, y)$ is essentially $1$ if $|x|, |y| \le R^{1/2}$, $0 \le t \le R$.
More precisely, we choose a large constant $C=C(n)$ so that $\Re F(t, x, y) \ge \frac{1}{2}$ in the above range. Also, as before, the Fourier transform (in space) of $F$ is essentially supported at frequencies $|\xi|, |\eta| \le R^{-1/2}$.

Pick roughly $R^{n-\frac 12}$ points $(x_k, y_k)$ which are spaced at distance $\sim R^{1/2}$ from each other on the surface $\{(x,y):|x|=|y|, \frac R2 \leq |x|\leq R\}$. Define $(\xi_k,\eta_k)=\frac 1R (x_k,y_k)$ so that $|\xi_k|^2-|\eta_k|^2=0$ and $|(\xi_k,\eta_k)|\sim 1$.

Take the following initial conditions
$$
\Gamma_0(x, y)=\sum e^{i (x \cdot \xi_k - y \cdot \eta_k)}F_0(x+x_k, y+y_k)
$$
so that the solution is
\begin{align*}
\Gamma(t, x, y)&=\sum e^{-i t \frac{ (|\xi_k|^2-|\eta_k|^2)}{2}}e^{i (x \cdot \xi_k - y \cdot \eta_k)}F(t, x+x_k- t \xi_k, y+y_k-t \eta_k)\\
&=\sum e^{i (x \cdot \xi_k - y \cdot \eta_k)}F(t, x+x_k- t \xi_k, y+y_k-t \eta_k).
\end{align*}

Since the $\sim R^{n-\frac 12}$ terms in $\Gamma_0$ are essentially orthogonal and each have $L^2$ norm $\sim R^{n/2}$, we get
$$
\|\Gamma_0\|_{L^2(dxdy)}\lesssim R^{n-\frac 14}.
$$
Moreover, since $|(\xi_k,\eta_k)|\sim 1$, there also holds
\begin{equation}\label{gahs}
   \|\Gamma_0\|_{H^s(dxdy)}\lesssim R^{n-\frac 14}.
\end{equation}

From the expression of $\Gamma$, we see that $$|\Gamma(t,x,y)|\gtrsim R^{n-\frac 12} \quad \text{ for }|(x, y) | \le \frac{1}{100}, R-R^{\frac{1}{2}}<t<R.$$
Therefore,
$$
\|\Gamma(t,x,x)\|_{L^p(dt)L^q(dx)}\gtrsim R^{\frac{1}{2p}} R^{n- \frac{1}{2}},
$$
so, recalling \eqref{gahs}, if
\begin{align*}
& \|  \Gamma(t, x, x) \|_{L^p(dt)L^q( dx )} \lesssim \| \Gamma_0(x, y)\|_{H^s(dx dy)},
\end{align*}
then $p \ge 2$.
From a similar argument to the one in subsection \ref{sec-lambdap} (i.e. only using $x_k$, $y_k$ for which all coordinates of
$\xi_k $ and $- \eta_k$ are $\ge \frac{1}{10n}$),
 $p\ge 2$ is also necessary for estimates of the form
$$
\||\nabla|_x^{\alpha} \Gamma(t, x, x) \|_{L^p(dt)L^q( dx )} \lesssim \| \Gamma_0(x, y)\|_{H^s(dx dy)}.
$$

\subsubsection{Necessity of $q \geq 2$.}

 Let $F(t, x, y)$ be the basic vertical tube solution of height $R$ (as \eqref{gammatube}). Let $m\gg1$. Choose roughly $R^{mn-\frac n2}$ points $x_k$ which are spaced at distance $\sim R^{\frac{1}{2}}$ in a large ball $B(0, R^m)$ of radius $R^m$ in $\mathbb R^n$. Fix a unit vector $\xi\in S^{n-1}$.

 We take initial conditions
 \begin{align*}
\Gamma_0( x, y)=e^{i x\cdot \xi}\sum F_0(x+x_k,y+x_k),
\end{align*}
so that the solution is
\begin{align*}
\Gamma(t, x, y)=e^{i x\cdot \xi}\sum F(t, x+x_k-t\xi, y+x_k).
\end{align*}
Note that $\Gamma(t,x,x)\gtrsim 1$ for $0\leq t\leq 1$ and $|x|\leq R^{m}$. Moreover, the Fourier transform of $\Gamma(t,x,x)$ is essentially supported in a $R^{-1/2}$ neighborhood of the point $\xi$ with $|\xi|=1$.

Then, the necessity of $q \ge 2$ follows from the same calculation as in subsection \ref{sec-lambdaq}.

\subsection{Proof of Theorem \ref{Gtheorem}}

The examples for $G$ are similar to those in previous subsections.

\subsubsection{Necessity of $p\geq 2$.}

First we take the basic ``vertical tube'' solution. Let
\begin{align*}
F_0(x, y, z)=  e^{-\frac{|x|^2+|y|^2+|z|^2}{2 C R}}
\end{align*}
so that
\begin{align}
e^{ \frac{it (\Delta_x +\Delta_y-\Delta_z)}{2}}F_0:&=
F(t, x, y, z) \notag\\&= \frac{1}{(1+\frac{it}{CR})^n(1-\frac{it}{CR})^{\frac n2}} e^{-\frac{|x|^2+|y|^2}{2(CR+it)}}e^{-\frac{|z|^2}{2(CR-it)}}. \label{Gtube}
\end{align}
The solution
$F(t, x, y,z)$ is essentially $1$ if $|(x,y,z)| \le R^{1/2}$, $0 \le t \le R$.
 Also, the Fourier transform (in space) of $F$ is essentially supported at frequencies $|(\xi,\eta,\zeta)| \le R^{-1/2}$.

Pick roughly $R^{\frac{3n-1}{2}}$ points $(x_k, y_k, z_k)$ which are spaced at distance $\sim R^{1/2}$ from each other on the surface $\{(x,y,z):|x|^2+|y|^2=|z|^2, \frac R2 \leq |x|,|y|\leq R\}$. Define $(\xi_k,\eta_k, \zeta_k)=\frac 1R (x_k,y_k,z_k)$ so that $$|\xi_k|^2+|\eta_k|^2=|\zeta_k|^2 \quad \text{ and }\quad  |(\xi_k,\eta_k,\zeta_k)|\sim 1.$$

Take the following initial conditions
$$
G_0(x, y,z)=\sum e^{i (x \cdot \xi_k + y \cdot \eta_k -z\cdot \zeta_k)}F_0(x+x_k, y+y_k, z+z_k)
$$
so that the solution is
\begin{align*}
&G(t, x, y,z)\\&=\sum e^{i (x \cdot \xi_k + y \cdot \eta_k -z\cdot \zeta_k)}F(t, x+x_k- t \xi_k, y+y_k-t \eta_k, z+z_k-t\zeta_k),
\end{align*}
since $|\xi_k|^2+|\eta_k|^2=|\zeta_k|^2$.

Since the $\sim R^{\frac{3n-1}{2}}$ terms in $G_0$ are essentially orthogonal and each has $L^2$ norm $\sim R^{3n/4}$, we get
$$
\|G_0\|_{L^2(dxdydz)}\lesssim R^{\frac{3n}{2}-\frac 14}.
$$
Moreover, since $|(\xi_k,\eta_k, \zeta_k)|\sim 1$, there also holds
\begin{equation}\label{Ghs}
   \|G_0\|_{H^s(dxdydz)}\lesssim R^{\frac{3n}{2}-\frac 14}.
\end{equation}

From the expression of $G$, we see that $$|G(t,x,y,z)|\gtrsim R^{\frac{3n-1}{2}} \quad \text{ for }|(x, y,z) | \le \frac{1}{100}, R-R^{\frac{1}{2}}<t<R.$$
Therefore,
$$
\|G(t,x,x,x)\|_{L^p(dt)L^q(dx)}\gtrsim R^{\frac{1}{2p}} R^{\frac{3n-1}{2}}.
$$
Recalling \eqref{Ghs}, if
\begin{align*}
& \| G(t, x, x, x) \|_{L^p(dt)L^q( dx )} \lesssim \| G_0(x, y,z)\|_{H^s(dx dydz)},
\end{align*}
then $p \ge 2$.
From a similar argument as in subsection \ref{sec-lambdap}, $p\ge 2$ is also necessary for estimates of the form
$$
\||\nabla|_x^{\alpha} G(t, x, x, x) \|_{L^p(dt)L^q( dx )} \lesssim \| G_0(x, y,z)\|_{H^s(dx dydz)}.
$$

\subsubsection{Necessity of $q \geq 2$.}

 Let $F(t, x, y, z)$ be the basic vertical tube solution of height $R$ (as \eqref{Gtube}). Let $m\gg1$. Choose roughly $R^{mn-\frac n2}$ points $x_k$ which are spaced at distance $\sim R^{\frac{1}{2}}$ in a large ball $B(0, R^m)$ of radius $R^m$ in $\mathbb R^n$. Fix a unit vector $\xi\in S^{n-1}$.

 We take initial conditions
 \begin{align*}
G_0( x, y,z)=e^{i (x+y-z)\cdot \xi}\sum F_0(x+x_k,y+x_k,z+x_k),
\end{align*}
so that the solution is
\begin{align*}
&G(t, x, y)\\&=e^{\frac{-it}{2}}e^{i (x+y-z)\cdot \xi}\sum F(t, x+x_k-t\xi, y+x_k-t\xi, z+x_k-t\xi).
\end{align*}

There are roughly $R^{mn - \frac{n}{2}}$ terms in the sum. The summands are essentially orthogonal and each term has $L^2$ norm
$\sim R^{3n/4}$, thus
\begin{align*}
\|G_0\|_{L^2(dxdydz)}\sim R^{\frac{n}{2}+\frac{mn}{2}}.
\end{align*}
On the other hand, each $F(t, x+x_k-t\xi, y+x_k-t\xi, z+x_k-t\xi)$ is essentially $1$ on a tube $T_k$ of radius $R^{1/2}$ and length $R$ in $3n+1$ dimensions, and rapidly decaying out of $T_k$. Note that at $t=0$, $T_k$ is centered at $(0, -x_k,-x_k,-x_k)$. Moreover, these tubes $T_k$ are in the same direction $(1,\xi,\xi,\xi)$ and hence disjoint. Therefore,
$|G(t,x,y,z)|\gtrsim 1$
on the union of the tubes $T_k$. In particular, $|G(t,x,x,x)|\gtrsim 1$ for $0\leq t\leq R$ and $x\in B(t\xi, R^m)$.
Thus
$$\|G(t, x, x,x)\|_{L^p([0, 1])L^q(dx)}\gtrsim  R^{\frac{mn}{q}}$$
(with a similar estimate for $|\nabla|^{\alpha}G(t, x, x,x)$),
while $\|G_0\|_{H^s(dxdy)}\sim R^{\frac{n}{2}+\frac{mn}{2}}$ and $m\gg 1$, so $q \ge 2$ is necessary.


\begin{thebibliography}{24}

\bibitem{C-P} T. Chen and  N. Pavlovi\'c, \emph{Derivation of the cubic NLS and Gross-Pitaevskii hierarchy from many body dynamics in  $d= 3$ based on spacetime norms},
Ann. H. Poincare, \textbf{15} (2014), 543--588.

\bibitem{C-H-P} T. Chen, Y. Hong and N. Pavlovi\'c,
\emph{Global Well-Posedness of the NLS System for Infinitely Many Fermions},
Archive for rational mechanics and analysis, April 2017, Volume \textbf{224}, Issue 1, pp 91--123.

\bibitem{C-H1} X. Chen and J. Holmer,
\emph{On the Klainerman-Machedon Conjecture of the Quantum BBGKY Hierarchy with Self-interaction}
, Journal of the European Mathematical Society, 2016 \textbf{ 18} , 1161-120.

\bibitem{C-H2} X. Chen and J. Holmer,
\emph{Correlation structures, Many-body Scattering Processes and the Deriva-
tion of the Gross-Pitaevskii Hierarchy}. Int. Math. Res. Not. 2016 (\textbf{10}), 3051--3110.

\bibitem{E-S-Y2} L. Erd\"os, B. Schlein and H.~T. Yau,
\emph{Derivation of the
   cubic non-linear Schr\"odinger equation from quantum dynamics of many-body systems}.
   Invent. Math. \textbf{167}, 515--614 (2007).



\bibitem{E-S-Y3}L. Erd\"os, B. Schlein and H.~T. Yau,
\emph{Derivation of the Gross-Pitaevskii equation for the dynamics of Bose-Einstein condensate}. Annals Math. \textbf{172}, 291--370 (2010).

\bibitem{FLLS} R. L. Frank, M. Lewin, E. H. Lieb, and R. Seiringer, \emph{Strichartz  inequality  for  orthonormalfunctions}, J. Eur. Math. Soc., (2014).

    \bibitem{FS}  R. L. Frank and J. Sabin, \emph{Restriction theorems for orthonormal functions, Strichartz inequalities and uniform Sobolev estimates}, American Journal of Mathematics
Johns Hopkins University Press
Volume \textbf{139}, Number 6, December 2017
pp. 1649--1691.



\bibitem{GM1} M. Grillakis and M. Machedon,  \emph{Pair excitations and the mean field approximation of interacting Bosons, II},  Communications in PDE, Vol \textbf{42}, No 1, 24--67 (2017).

\bibitem{GM2} M. Grillakis and M. Machedon, \emph{Uniform in $N$ estimates for a Bosonic system of Hartree-Fock-Bogoliubov type},  Communications in PDE,
Volume \textbf{44}, Number 12,  2019, pp. 1431--1465.

\bibitem{K-MMM}
 S. Klainerman and M. Machedon, \emph{ On the uniqueness of solutions to the
    Gross-Pitaevskii
hierarchy}.  Comm. Math.\ Phys.\ \textbf{279}, 169--185 (2008).

\end{thebibliography}
\end{document}